\documentclass{amsart}
\usepackage{amssymb,amsfonts,amsthm,amsmath,amscd}
\numberwithin{equation}{section}

\newtheorem{theorem}{Theorem}[section]
\newtheorem{lemma}[theorem]{Lemma}
\newtheorem{corollary}[theorem]{Corollary}
\newtheorem{proposition}[theorem]{Proposition}
\theoremstyle{remark}
\newtheorem{remark}{Remark}[section]

\def\S{{\mathbb{S}}}

\def\eps{\epsilon}

\def\const{\qopname \relax o{const}}
\def\sub{\qopname \relax o{sub}}

\def\ra{\rightarrow}
\def\tr{\qopname \relax o{tr}}

\def\R{{\mathbb{R}}}
\def\N{{\mathbb{N}}}
\def\Z{{\mathbb{Z}}}

\def\C{{\mathbb{C}}}

\def\min{\qopname \relax o{min}}
\def\diag{\qopname \relax o{diag}}

\def\tr{\qopname \relax o{Tr}}

\def\Res{\qopname \relax o{Res}}
\def\det{\qopname \relax o{det}}

\def\ls{\left[}
\def\rs{\right]}

\def\spec{\qopname \relax o{spec}}

\begin{document}

\title{Lower Order Terms in Szeg\"o Theorems
on Zoll Manifolds}

\author{Dimitri Gioev}
\address{Department of Mathematics, University of Pennsylvania, 
DRL 209 South 33rd Street, Philadelphia, PA 19104--6395}
\curraddr{Department of Mathematics, Courant Institute
of Mathematical Sciences,
New York University, 251 Mercer Street, New York, NY 10012--1185}
\email{gioev@math.upenn.edu, gioev@cims.nyu.edu}
\thanks{The author gratefully acknowledges a full support for the academic year 2002--03 from the Swedish Foundation for International Cooperation in Research and Higher Education (STINT), Dnr.~PD2001--128.}

\date{May 25, 2002}



\begin{abstract}
We give an outline of the computation of the third order term in a
generalization of the Strong Szeg\"o Limit Theorem
for a zeroth order pseudodifferential operator (PsDO) on a Zoll manifold
of an arbitrary dimension, see \cite{G1} for the detailed proof. 
This is a refinement of a result 
by V.~Guillemin and K.~Okikiolu
who have computed the second order term in \cite{GO}.
An important role in our proof 
is played by a certain combinatorial identity which generalizes the formula of 
G.~A.~Hunt and F.~J.~Dyson to an arbitrary natural power, see \cite{G2}.
This identity is a different form of 
the renowned Bohnenblust--Spitzer combinatorial theorem
which is related to the maximum of a random walk with i.i.d. steps on the real line.

A corollary of our main result is a fourth order Szeg\"o type
asymptotics for a zeroth order PsDO on the unit circle, which
in matrix terms gives
a fourth order asymptotic formula
for the determinant of the truncated 
sum $P_n(T_1+T_2D)P_n$
of a Toeplitz matrix $T_1$ with the product 
of another Toeplitz matrix $T_2$ 
and a diagonal matrix $D$
of the form $\diag(\cdots,\frac13,\frac12,1,0,1,\frac12,\frac13,\cdots)$. 
Here $P_n=\diag(\cdots,0,1,\cdots,1,0,\cdots)$, $(2n+1)$ ones.
\end{abstract}

\maketitle

\section{Introduction and main results}
The main motivation for this work was to find an
explicit formula for a ``Szeg\"o--regularized'' determinant
of a zeroth order pseudodifferential operator (PsDO)
on a Zoll manifold introduced in \cite{GOs,O2}, see Remark~\ref{rem_appl}.
Our main result, Theorem~\ref{thm_main}, is valid for any dimension
$d\in\N$. In the case $d=2$, Theorem~\ref{thm_main}
gives such a formula.

In this paper we find a third order
generalization of the Strong Szeg\"o Limit Theorem (SSLT)
 for a zeroth order PsDO on the unit circle 
(Theorem~\ref{thm_onedim}), 
and on 
a Zoll manifold of an arbitrary dimension
$d\in\N$ (Theorem~\ref{thm_main}).
We give also an outline of the proof of Theorem~\ref{thm_main}
and derive Theorem~\ref{thm_onedim} from the former. 
The detailed proof of Theorem~\ref{thm_main}
can be found in 
\cite[Chapter~1]{G}
and will be published in a forthcoming paper \cite{G1}.
In Section~3, we sketch a proof of the 
generalized Hunt--Dyson combinatorial
formula (Theorem~\ref{thm_comb}),
see \cite{G2} for the detailed proof. 
We give also a brief review of the related work.

Recall that $M$ is called a {\em Zoll manifold\ }if it is compact, closed
and such that the geodesic flow
on $M$ is simply periodic with period $2\pi$. 
The unit circle and the standard sphere of any dimension
are Zoll manifolds.
A second order generalization of the SSLT for a Zoll manifold
of any dimension has been obtained by V.~Guillemin
and K.~Okikiolu \cite{GOs,GO}, 
see also an important preceding work \cite{O} by K.~Okikiolu
for the case of the two- and three-dimensional sphere.
The proofs in \cite{O,GOs,GO} use a combinatorial
identity due to G.~A.~Hunt and F.~J.~Dyson and proceed in
the spirit of the combinatorial proof of the classical SSLT 
by M.~Kac \cite{K}. See also \cite{GO3,O2} 
where the combinatorial approach and the Hunt--Dyson
formula are used in a different setting
to obtain
a second order generalization of the SSLT 
for a manifold with the set of closed geodesics
of measure zero in the unit cotangent bundle. 

In the proof of Theorem~\ref{thm_main} we use the method of \cite{GO}.
A central role in our proof is
played by a certain combinatorial identity
which generalizes the Hunt--Dyson 
formula mentioned above to an arbitrary natural power.
We call this identity the {\em generalized Hunt--Dyson formula (gHD)}, 
see Theorem~\ref{thm_comb}.

After having discovered and proved the formula gHD
we realized that it is related to 
another combinatorial theorem, which has a long history. 
The mentioned theorem 
is a result due to H.~F.~Bohnenblust
that appeared in an article by F.~Spitzer
on random walks \cite[Theorem~2.2]{S1}, 
and is now commonly known as the
{\em Bohnenblust--Spitzer theorem (BSt)}.
The characteristic function of the maximum of a random walk
with independent identically distributed steps has been computed
in \cite{S1} with the help of the BSt. 
On the other hand, using the gHD, we can compute \cite{G2}
the moment of an arbitrary order of
the mentioned maximum (note that 
the usual Hunt--Dyson formula allows one to compute
only the expected value of the maximum, see \cite{K}).
This indicates that the gHD and the BSt should be essentially the same. 
And indeed this is the case: 
it turns out that the gHD can be derived from the BSt,
and vice versa, see \cite{G2}. 

We have found the formula gHD being unaware of the BSt.
In \cite[Chapter~2]{G}, a proof of the 
gHD ``from scratch'' can be found. 
This together with a 
 a derivation of the BSt from
the the gHD in \cite{G,G2}
gives a new proof of the BSt.
There are known
at least four other proofs: the original one \cite{S1,F}, 
a proof by G.~Baxter \cite{B1,B2},
by J.~G.~Wendel \cite{We}, and finally, a unifying approach of G.-C.~Rota
in the framework of (Glen) Baxter algebras \cite{Ro,RoSm}.
See Section~3 for 
more combinatorial details.

Let $\det$ denote the determinant of a finite rank operator.
In Section~2, we give explicit asymptotic formulas for
$\log\det{}P_nBP_n$, $n\ra\infty$, 
where $P_n$ is the projection onto the first $n$ 
eigenspaces of the Laplace--Beltrami
operator on a Zoll manifold, and $B$ is an arbitrary zeroth order PsDO,
see Corollary~\ref{cor4} (for $d=1$) and Corollary~\ref{cor5} ($d\geq2$).

Let $\S^1$ be the unit circle $\R/2\pi\Z$. 
Denote by $P_n$, $n\in\N$, the orthogonal
projection from $L^2(\S^1)$ 
to the subspace spanned by $\{e^{ikx}\}_{|k|\leq{}n}$. 
For a function $f\in{}L^1(\S^1)$ 
denote its $k$th
Fourier coefficient by
\begin{equation}
\label{eq_Fcf}
            \widehat{f}_k 
  := \int_0^{2\pi} f(x) e^{-ikx}\,\frac{dx}{2\pi},\qquad k\in\Z.
\end{equation}
Let $b(x)$ be a strictly positive function on $\S^1$
such that $\sum_{k\in\Z}|k|\,|\widehat{(\log{b})}_k|^2<\infty$.
Denote by $B$ the operator of multiplication by
$b$ acting in $L^2(\S^1)$.
Recall that the matrix representation of the operator $B$
is the Toeplitz matrix $(\widehat{b}_{j-k})_{j,k\in\Z}$.
The classical Strong Szeg\"o Limit Theorem (SSLT)~\cite{Sz}
states that
$$
   \tr\log P_nBP_n = \tr P_n(\log{B})P_n
              +\sum_{k=1}^\infty k\,\widehat{(\log b)}_k\,
                       \widehat{(\log b)}_{-k}
              +o(1),\qquad n\ra\infty.
$$
Here $\tr\log P_nBP_n=\log\det P_nBP_n$
and $\tr P_n(\log{B})P_n=(2n+1)\int_0^{2\pi} \log b(x)\,\frac{dx}{2\pi}$.
It has been shown by H.~Widom that the remainder
is $O(n^{-\infty})$ if $b(x)\in{}C^{\infty}(\S^1)$, see \cite{Wajm}.

Let $M$ be a closed manifold of dimension $d\in\N$. 
Let $\Psi^m(M)$, $m\in\Z$, denote
the space of the classical PsDO of order $m$ on $M$. 
Recall that for a given $G\in\Psi^m(M)$,
its principal symbol $\sigma_m(G)$ and subprincipal symbol $\sub(G)$
are well-defined on $T^*M$.
Let $S^*M$ be the unit cotangent bundle and denote
by $dxd\xi$ the standard measure on $S^*M$ divided by $(2\pi)^d$.

The simplest form of our result is for $d=1$ and the  
function $\psi(u)=\log{}u$ which is analytic in
a neighborhood of $u=1$.
\begin{theorem}
\label{thm_onedim}
 Let $M=\S^1$ and $P_n$ be the projection
on the linear span of $\{e^{ikx}\}_{|k|\leq{}n}$.
Let $B\in\Psi^0(M)$ and assume that $\sigma_0(B)$
is strictly positive and 
a certain symbolic norm of $I-B$
is sufficiently small. Then $\log{}B\in\Psi^0(M)$ and 
the following holds
\begin{equation}
\label{eq_f}
\begin{aligned}
    \tr\log\,&P_nBP_n = \tr P_n(\log B)P_n 
                + \frac12\int_{S^*M}\sum_{k=1}^\infty k\,\widehat{(\sigma_0(\log{B}))}_k
                   \,\widehat{(\sigma_0(\log{B}))}_{-k}\,dxd\xi\\
              &+\frac1n\cdot
                     \frac12\int_{S^*M}\sum_{k=1}^\infty k\,\widehat{(\sigma_0(\log{B}))}_k
                   \,\widehat{(\sub(\log B))}_{-k}\,dxd\xi
         +O\bigg(\frac1{n^2}\bigg),\quad n\ra\infty.
\end{aligned}
\end{equation}
\end{theorem}
In \eqref{eq_f} the argument $(x,\xi)\in{}S^*M$ is omitted 
for brevity. By the Fourier coefficient we mean the
following: for a fixed $(x,\xi)\in{}S^*M$ compute the Fourier
coefficient in accordance with \eqref{eq_Fcf} 
along the unit circle being
the closed geodesic starting at $(x,\xi)$.

Let us denote $b_0:=\sigma_0(B)$ and $b_{{\rm sub}}:=\sub(B)$.
Recall that for an analytic $f$, 
$\sigma_0(f(B))=f(\sigma_0(B))$ and $\sub f(B)=f^\prime(\sigma_0(B))\sub(B)$. 
Then \eqref{eq_f} takes the form
\begin{equation}
\label{eq_f1}
\begin{aligned}
    \tr\log P_nBP_n &= \tr P_n(\log B)P_n 
                + \frac12\int_{S^*M}\sum_{k=1}^\infty k\,\widehat{(\log{b_0})}_k
                   \,\widehat{(\log{b_0})}_{-k}\,dxd\xi\\
              &+\frac1n\cdot
                    \frac12\int_{S^*M}\sum_{k=1}^\infty k\,\widehat{(\log{b_0})}_k
                   \,\widehat{(b_{{\rm sub}}/b_0)}_{-k}\,dxd\xi
         +O\bigg(\frac1{n^2}\bigg),\quad n\ra\infty.
\end{aligned}
\end{equation}
\begin{remark}
Theorem~\ref{thm_onedim} is a particular case 
of Theorem~\ref{thm_main} below, see Remark~\ref{rem_pfdim1}.
\end{remark}

Let us fix some notations and then state the result for the 
higher dimensional case.
Let $M$ be a Zoll manifold of dimension $d\in\N$.
Let $\Delta$ denote the Laplace--Beltrami operator on $M$.
It is known \cite{DG} that there exists a constant $\alpha\in\R$
such that the spectrum of  $\sqrt{-\Delta}$
lies in bands around the points $k+\frac{\alpha}{4}$, $k\in\N$.
Moreover, in has been shown in \cite{CdV}
that there exists $A_{-1}\in\Psi^{-1}(M)$
such that $[\Delta,A_{-1}]=0$
and the spectrum of the operator 
\begin{equation}
\label{eq_A}
     A:= \sqrt{-\Delta}-\frac{\alpha}4 - A_{-1}
\end{equation}
is $\N$.
Let $P_n$, $n\in\N$, denote the projection from $L^2(M)$
onto the subspace spanned by the eigenfunctions
of $A$ corresponding to the eigenvalues $1,2,\cdots,n$. 
Let $dxd\xi$ be the defined above measure on 
$S^*M:=\{(x,\xi):\sigma_1(A)(x,\xi)=1\}$. 
Following \cite{GOs,GO}, for a $B\in\Psi^0(M)$ and $t\in\R$ 
introduce the operator
\begin{equation}
\label{eq_Bt}
       B^t:=e^{-itA}Be^{itA}.
\end{equation}
By Egorov's theorem, $B^t\in\Psi^0(M)$, and also 
$$
     \sigma_0(B^t)(x,\xi)=\sigma_0(B)(\Theta^t(x,\xi)),
$$ 
where $\Theta^t$ stands for the shift by $t$ units 
along geodesic flow.
Note that because $\sub(A)=\const$ the following also holds \cite[Lemma~2.2]{Gui}
$$
     \sub(B^t)(x,\xi)=\sub(B)(\Theta^t(x,\xi)).
$$ 
Because $\spec(A)=\N$, the operator $B^t$ is periodic
in $t$ with period $2\pi$.

The result we are about to state, and also Theorem~\ref{thm_onedim} above,
gives information on the asymptotic behavior of 
\begin{equation}
\label{eq_trdiff}
      \tr \psi(P_nBP_n) - \tr P_n\psi(B)P_n,\qquad n\ra\infty,
\end{equation}
for the analytic in a neighborhood of $1$ function $\psi(u)=\log{}u$.
In \cite{G,G1} we obtain the corresponding statements for
an arbitrary analytic function $\psi(u)$. These results
involve certain maps $W$, $\Phi$, $\Upsilon$
whose action on an analytic function is a continuous function
from $\C^j$ to $\C$ for some $j\in\N$.
The map 
$$
      W[\psi](y_1,y_2):= \frac12\int_{y_1}^1\int_{y_2}^1
      \frac{\psi^\prime(u_1)-{\psi^\prime(u_2)}}
{u_1-u_2}\,du_1du_2,\qquad y_1,y_2>0,
$$ 
was obtained by A.~Laptev, D.~Robert
and Yu.~Safarov in \cite{LRS} (here $j=2$). 
Earlier, an equivalent to $W$ map 
has been obtained 
by H.~Widom \cite{W0,Wjfa}
in a two-term Szeg\"o type asymptotics 
for integral operators with discontinuous symbols.
It has been noticed by H.~Widom, 
and also by the authors \cite{LRS,O2},
that for the function $\psi(u)=\log{}u$, 
the action $W[\log]$ takes a simple form
$$
   W[\log](y_1,y_2)= -\frac12\log y_1 \log y_2.
$$
This explains the fact that the 
map $W$ is not required in the next theorem. 

Now define for $y_1,y_2,y_3<1$,
\begin{equation}
\label{eq_Phi}
\begin{aligned}
     \Phi[\log](y_1,y_2,y_3)&= y_3\int_0^{y_1}\!\!\int_0^{y_2}\bigg[
      \frac{\log(1-u_1)}{u_1(u_1-u_2)(u_1-u_3)}\\
     &-\frac{\log(1-u_2)}{u_2(u_1-u_2)(u_2-u_3)}
    +\frac{\log(1-y_3)}{y_3(u_1-y_3)(u_2-y_3)}\bigg]
\,du_1\,du_2.
\end{aligned}   
\end{equation}
For the map $\Phi$, $j=3$.
The expression for $\Upsilon[\log]$ 
appearing in the statement of the next theorem
is complicated and is not given here, see \cite[(1.1.31)]{G} and \cite{G1}
for details.
It depends only on the principal symbol $\sigma_0(B)$
and involves 
certain Poisson brackets of the type $\{\sigma_0(B^t),\sigma_0(B^s)\}$,
$0\leq{}t,s\leq2\pi$.

For a function $f\in{}C^\infty(S^*M)$ introduce
the $k$th Fourier coefficient along the closed geodesic
of length $2\pi$ starting at a given point $(x,\xi)$
$$
        \widehat{f}_k(x,\xi):= \int_0^{2\pi} e^{-ikt} f(\Theta^t(x,\xi))\,
     \frac{dt}{2\pi},\qquad k\in\Z.
$$
For an arbitrary $B\in\Psi^0(M)$, let us write $b_0:=\sigma_0(B)$, 
$b_{{\rm sub}}:=\sub(B)$, $b_0^t(x,\xi):=b_0(\Theta^t(x,\xi))$,
and omit the argument $(x,\xi)\in{}S^*M$ for brevity.
\begin{theorem}
\label{thm_main}
Let $M$ be a Zoll manifold of dimension $d\in\N$. 
Let $A$ be defined by \eqref{eq_A}.
Assume that $\sigma_1(A)(x,\xi)=\sigma_1(A)(x,-\xi)$ for all $(x,\xi)\in{}T^*M$.
Let $B\in\Psi^0(M)$ be such that
$\sigma_0(B)$ is strictly positive and a certain symbolic norm of $I-B$
is sufficiently small. Then $\log{B}\in\Psi^0(M)$ and 
the following holds
\begin{equation}
\label{eq_th2}
\begin{aligned}
    \tr\log\,&P_nBP_n = \tr P_n(\log B)P_n \\
                &+n^{d-1}\cdot\frac12\int_{S^*M}
\,\sum_{k=1}^\infty k\,\widehat{(\log b_0)}_k
                   \,\widehat{(\log b_0)}_{-k}\,dxd\xi\\
              &+n^{d-2}\cdot\frac12\int_{S^*M}\bigg(
                    \sum_{k=1}^\infty k\,\widehat{(\log b_0)}_k
                   \,\widehat{(b_{{\rm sub}}/b_0)}_{-k}\\
              &\qquad\,\,+ (d-1)
\bigg[
\sum_{k=1}^\infty\Big(k^2+(1+\alpha/2) k\Big)\,
      \widehat{(\log b_0)}_k
                   \,\widehat{(\log b_0)}_{-k}\\
               &\qquad\qquad
\,+\sum_{k,l=1}^\infty  k\,l\,\int_0^{2\pi}\!\!
\int_0^{2\pi}\!\!\int_0^{2\pi} e^{ik(t-r)+il(s-r)}
    \Phi[\log](b_0^t,b_0^s,b_0^r)\,\frac{dt}{2\pi}
\frac{ds}{2\pi}\frac{dr}{2\pi}\bigg]\\
 &\qquad\,\,+ \Upsilon[\log](b_0)\bigg)\,dxd\xi\\
             &+O(n^{d-3}),\qquad n\ra\infty.
\end{aligned}
\end{equation}
\end{theorem}
\begin{remark}
The {\em existence\ }of a full expansion of the type \eqref{eq_th2}
has been proven in \cite{GOs}. Explicit formulas
for the coefficients of the first first two terms
in \eqref{eq_th2} (and also of the $\log{}n$ term which is a part
of $\tr P_n(\log{}B)P_n$) 
have been obtained in \cite{GOs,GO}.
\end{remark}
\begin{remark}
The formula \eqref{eq_Phi} for $\Phi[\log]$ 
has a structure similar to the
coefficient in the third asymptotic term in a Szeg\"o
type expansion for convolution operators obtained by
R.~Roccaforte in \cite{R}. A version of the Bohnenblust--Spitzer 
combinatorial theorem is used in the proof
in \cite{R}. A second order Szeg\"o type expansion for convolution operators
was established by H.~Widom in \cite{Wn} 
with the help of the usual Hunt--Dyson combinatorial
formula \eqref{CH2_s2e10orig}. A full asymptotic expansion
for convolution operators has been found in \cite{Wbook}.
\end{remark}
\begin{remark}
\label{rem_pfdim1} 
For $d=1$, the square bracket in \eqref{eq_th2} disappears.
Also $\Upsilon$ vanishes, because all the Poisson brackets
vanish in this case (for each of the two cotangent directions
the angle does not change and
$\sigma_0(B)$ is homogeneous of degree $0$ in $\xi$). Theorem~\ref{thm_onedim} 
follows.
\end{remark}
\begin{remark}
As we have noticed above,
Theorem~\ref{thm_onedim} and Theorem~\ref{thm_main} 
give lower order corrections 
to \eqref{eq_trdiff} for $\psi(u)=\log{}u$, as $n\ra\infty$. 
In the next section we find explicit formulas 
for the coefficients in the asymptotic 
expansion of $\tr\log{}P_nBP_n$, as $n\ra\infty$.
These expansions are fourth order for $d=1,2,3$
and third order for $d\geq4$, as $n\ra\infty$.  
\end{remark}
\begin{proof}[Proof of Theorem~\ref{thm_main}: an outline]
We start by expanding $\psi(u)=\log{}u$ 
in a power series about the point $1$.
After that it suffices to prove \eqref{eq_th2} for $\psi(u)=u^m$
for all $m\in\N$. 

Let us recall the method of \cite{GO}. 
Let $\pi_k$, $k\in\N$, be the projection on the $k$th eigenspace
of the operator $A$ and set $\pi_k:=0$ for $k\leq0$.
Then $P_n=\sum_{k=1}^n\pi_k$ for $n\in\N$, and we set
$P_n:=0$, $n\leq0$.
Recall that because $M$ is a Zoll manifold, for an arbitrary
$B\in\Psi^0(M)$, the operator $B^t$, $t\in\R$, defined by \eqref{eq_Bt}
is $2\pi$-periodic.
Therefore one can introduce the Fourier expansion
$
          B = \sum_{j\in\Z} B_{j}
$
where $B_{j}\in\Psi^0(M)$, $j\in\Z$, is defined by
$$
\begin{aligned}
           B_j &= \sum_{k=1}^\infty 
                                     \pi_{k+j}\,B\,\pi_k 
                 &=  \frac1{2\pi}\int_0^{2\pi} e^{ij{}t}\,
                                  e^{-itA}Be^{itA}\,dt.
\end{aligned}
$$
For $m\in\N$ and $j_1,\cdots,j_m\in\Z$
introduce the notation
$$
  M_m(\overline{j}):= \min(0,j_1,j_1+j_2,\cdots,j_1+\cdots+j_m).
$$
The following commutation relation is of central importance here
\begin{equation}
\label{eq_cr}
      B_j P_n = P_{n+j} B_j,\qquad n\in\N,\,\, j\in\Z.
\end{equation}
Using \eqref{eq_cr} one moves all the projectors to the left in the expression
$$
   (P_nBP_n)^m = \sum_{j_1,\cdots,j_m}P_nB_{j_1}P_n
             B_{j_2}P_n\cdots P_nB_{j_m}P_n
$$ 
obtaining $P_n B^m P_n$ plus a correction. This implies for all $n\in\N$
\begin{equation}
\label{CH1_s5e4}
\begin{aligned}
    \tr\,(P_n B &P_n)^m - \tr P_n B^m P_n\\
       &= -\sum_{j_1+\cdots+j_m=0}\,\tr\,\big(
          (P_n-P_nP_{n+j_1}\cdots{}P_{n+j_1+\cdots+j_m}\big)
                B_{j_1}\cdots B_{j_m}\big)\\
&=-\sum_{j_1+\cdots+j_m=0}\,\,
          \sum_{k=0}^{M_m({\overline{j}})+1}
               \tr(\pi_{n+k}\,B_{j_1}\cdots B_{j_m}).
\end{aligned}
\end{equation}
Next, it is important that 
for any $G\in\Psi^0(M)$, $M$ being a Zoll manifold,
there exists a full asymptotic expansion for $\tr(\pi_k{}G)$, as $k\ra\infty$,
see Lemma~\ref{lem_CdV} below.
This result is due to Y.~Colin de Verdi\`ere \cite{CdV}.
The coefficients in that expansion
are certain Guillemin--Wodzicki residues.
Recall that for any compact closed manifold $M$
of dimension $d\in\N$
the Guillemin--Wodzicki residue 
of a pseudodifferential operator $G\in\Psi^m(M)$ of order $m\in\Z$
is defined by
$$
   \Res(G):= \int_{S^*M} \sigma_{-d}(G)(x,\xi)\,dxd\xi.
$$
For an arbitrary $G\in\Psi^0(M)$ denote
$$
      R_l(G) :=\Res (A^{-d+l}G),\qquad l=0,1,2,\cdots.
$$
\begin{lemma}
\label{lem_CdV}
Let $M$ be a Zoll manifold of dimension $d\in\N$.
Assume $G\in\Psi^0(M)$.
Then for any $N=0,1,2,\cdots$, 
there exists $C_N<\infty$ such that
\begin{equation}
\label{eq_CdV}
         \Big| \tr(\pi_kG)
   -\sum_{l=0}^{N} k^{d-1-l}\,R_l(G) \Big|\leq C_N\,k^{d-2-N}, \quad k\in\N.
\end{equation}
\end{lemma}
\noindent
See \cite{CdV} and \cite[Appendix]{GO} for the proof of \eqref{eq_CdV}.

Now assume that for some $K<\infty$ we
have $B_k=0$ for $|k|>K$ (this is not assumed in the full proof in \cite{G1}).
Then \eqref{CH1_s5e4} and Lemma~\ref{lem_CdV} imply
\begin{equation}
\label{eq_big}
\begin{aligned}
    \tr\,(&P_n B P_n)^m - \tr P_n B^m P_n \\
       &= \sum_{j_1+\cdots+j_m=0}
             \bigg\{  n^{d-1}\cdot{}
  M_m({\overline{j}})\,\int_{S^*M}\sigma_0(B_{j_1}\cdots B_{j_m})\,dxd\xi \\
&\qquad+n^{d-2}\cdot\bigg[\frac{d-1}{2}\Big(
   \big[M_m({\overline{j}})\big]^2+\big(1+\frac{\alpha}{2}\big)\,M_m({\overline{j}})\Big)\,
     \int_{S^*M}\sigma_0(B_{j_1}\cdots B_{j_m})\,dxd\xi \\
&\qquad\qquad\qquad\qquad\quad\qquad
+M_m({\overline{j}})\,\int_{S^*M}\sub(B_{j_1}\cdots B_{j_m})\,dxd\xi
\bigg]\bigg\}\\
&\qquad+O(n^{d-3}),\qquad n\ra\infty.
\end{aligned}
\end{equation}
Now the symmetrization argument from \cite{K,S1,O,GO,GO3,O2} comes into play.
Note that the domain of summation over $j_1,\cdots,j_m$
in \eqref{eq_big} is symmetric with respect to the permutations of the $j$'s.
Also the factors involving 
$$
      \sigma_0(B_{j_1}\cdots{}B_{j_m}) = \sigma_0(B_{j_1})\cdot\ldots\cdot\sigma_0(B_{j_m}) 
$$
are symmetric. Therefore one can write
$$
\begin{aligned}
      \sum_{j_1+\cdots+j_m=0}
            &M_m({\overline{j}})\int_{S^*M}\sigma_0(B_{j_1}\cdots B_{j_m})\,dxd\xi\\
    &= \sum_{j_1+\cdots+j_m=0}
            \bigg(\frac1{m!}\sum_{\tau\in{}S_m}M_m({\overline{j}_\tau})\bigg)
   \int_{S^*M}\sigma_0(B_{j_1}\cdots B_{j_m})\,dxd\xi,  
\end{aligned}
$$
where $S_m$ is the set of all permutations of $j_1,\cdots,j_m$
and ${\overline{j}_\tau}:=(j_{\tau_1},\cdots,j_{\tau_m})$.
After an application of the usual Hunt--Dyson combinatorial formula
(see \eqref{CH2_s2e10orig} in Section~3)
$$
        \sum_{\tau\in{}S_m}M_m({\overline{j}_\tau})
=\sum_{\tau\in{}S_m}\sum_{k=1}^m
\frac{\min(0,j_{\tau_1}+\cdots+j_{\tau_k})}{k},
$$
all summations but one in $\sum_{j_1+\cdots+j_m=0}$
become free and can be carried out.

Now we also see the two difficulties with the computation
of the $n^{d-2}$ term. The first difficulty is in computing 
the minimum raised to power $2$. Here one needs a
generalization of the Hunt--Dyson formula, see Theorem~\ref{thm_comb}.
The second difficulty is that the factor
$$
     \int_{S^*M}\sub(B_{j_1}\cdots{}B_{j_m})\,dxd\xi
$$
is generally speaking {\em not\ }symmetric with respect to the
permutations $j_1,\cdots,j_m$, and the straightforward
symmetrization as above fails.
Indeed, recall that 
\begin{equation}
\label{CH1_s9e10}
\begin{aligned}
  \sub(B_{j_1}\cdots{}B_{j_m}) &= \sum_{k=1}^m\sub(B_{j_k})
                       \prod_{\genfrac{}{}{0pt}{}{p=1}{p\neq{}k}}^m\sigma_0(B_{j_p}) \\
                       &+ \frac1{2i}\,\sum_{1\leq{}k<l\leq{}m}
                               \big\{\sigma_0(B_{j_k}),\sigma_0(B_{j_l})\big\}\,
           \prod_{\genfrac{}{}{0pt}{}{p=1}{p\neq{}k,p\neq{}l}}^m\sigma_0(B_{j_p}).
\end{aligned}
\end{equation}
where $\{\cdot,\cdot\}$ stands for the Poisson bracket.
We see that the first sum 
in \eqref{CH1_s9e10} is symmetric, whereas
the second one is generally speaking not (not even after 
the integration over $S^*M$). However, each of the $m(m-1)/2$
terms in the second sum on the right-hand side of \eqref{CH1_s9e10}
possesses a {\em partial\ }symmetry. More precisely, we 
are allowed to permute the principal symbols 
that do not enter the Poisson bracket. 
This fact can be utilized in a modification of the symmetrization
procedure. The corresponding contribution 
is lengthy, this is where the functional $\Upsilon$ arises.
\end{proof}
\begin{remark}
We would like to mention here that in the computation of $\Upsilon$ 
one really needs the information on the {\em set of values\ }
$\big\{M_m({\overline{j}_\tau})\big\}_{\tau\in{}S_m}$ counted with 
multiplicities, 
that is the Bohnenblust--Spitzer theorem (see Section~3 for its statement),
and not just a formula for $\sum_{\tau\in{}S_m}
[M_m({\overline{j}_\tau})]^p$ for $p=1,2$.
\end{remark}
\section{Explicit asymptotic formulas for $\log\det{}P_nBP_n$, as $n\ra\infty$.}
Theorem~\ref{thm_onedim} and \ref{thm_main} give an expression
for $\tr\,\log{P_nBP_n}=\log\det{}P_nBP_n$ as a sum of $\tr{}P_n(\log{B})P_n$
and the two lower order corrections, as $n\ra\infty$.
We would like to compute the coefficients
in the asymptotic expansion of $\log\det{}P_nBP_n$,
as $n\ra\infty$.

An auxiliary asymptotic expansion for 
\begin{equation}
\label{eqstrace}
   \tr P_nGP_n=\sum_{k=1}^n\tr(\pi_kG),\qquad n\ra\infty,
\end{equation}
where $G\in\Psi^0(M)$, is given in Proposition~\ref{prop3} below.
A corresponding result for the coefficients of $n^d$, $n^{d-1}$
and $\log{}n$ for $d\geq2$ can be found in
\cite[after Lemma~0.2]{GO}.
To prove Proposition~\ref{prop3}, 
we sum over $k=1,\cdots,n$ in \eqref{eq_CdV}, see \cite{G1} for details.
The subtle point is the constant coefficient in \eqref{eqstrace}, 
which we
need for dimension $d=1,2$. The terms of all orders
in \eqref{eq_CdV}, and also the possible rapidly decaying term, 
will contribute to it. 
Let us therefore for $d=1,2$ 
make an additional assumption
\begin{equation}
\label{eq_assu}
   \sum_{l=0}^\infty{}\big|R_l(G)\big|<\infty
\end{equation}
and set, for all $k\in\N$,
\begin{equation}
\label{eq_eps}
   \eps_k(G):= \tr(\pi_kG)
   -\sum_{l=0}^{+\infty} k^{d-1-l}\,R_l(G).
\end{equation}
If \eqref{eq_assu} holds, the series in \eqref{eq_eps} (and also in \eqref{eq_C}
below)
is absolutely convergent, and also for any $N\in\N$ there exists $c_N<\infty$
such that $|\eps_k|\leq{}c_Nk^{-N}$, $k\in\N$. 
Set
\begin{equation}
\label{eq_C}
   C(G):= \sum_{k=1}^\infty \eps_k(G).
\end{equation}
Note that in \eqref{eq_C} and also in the 
proposition below there only appears 
the series $\sum_{l=0}^\infty{}R_l(G)$.
However we need the absolute convergence \eqref{eq_assu}
in the proof of the remainder estimate.
 
Let $\gamma$ denote the Euler constant and $\zeta$ the Riemann
zeta function.
\begin{proposition}
\label{prop3}
Let $M$ be a Zoll manifold of dimension $d\in\N$.
Let $P_n$ be as above and
assume that $G\in\Psi^0(M)$. For $d=1,2$, assume in addition
that \eqref{eq_assu} holds, and let $C(G)$ be defined by \eqref{eq_C}. 
Then the following holds, as $n\ra\infty$,\newline
(i) for $d=1$,
$$
\begin{aligned}
 \tr\,&P_nGP_n = n\cdot R_0(G) + \log n\cdot R_1(G)\\
 &+ 
  \Big(C(G)+\gamma\,R_1(G)+\sum_{l=2}^\infty\zeta(l)\,R_{l}(G)\Big)
+\frac1n\cdot\Big(\frac12\,R_1(G) - R_2(G)\Big)  + O\Big(\frac1{n^{2}}\Big),
\end{aligned}
$$
(ii) for $d=2$,
$$
\begin{aligned}
 \tr P_nGP_n &= n^2\cdot\frac12\,R_0(G) + 
  n\cdot\Big(\frac12\,R_0(G) + R_1(G)\Big)
+ \log n\cdot R_2(G)\\
 &+ 
  \Big(C(G)+\gamma\,R_2(G)+\sum_{l=2}^\infty\zeta(l)\,R_{l+1}(G)\Big)
  + O\Big(\frac1n\Big),
\end{aligned}
$$
(iii) for $d\geq3$,
$$
\begin{aligned}
 \tr\,&P_nGP_n = n^d\cdot\frac1d\,R_0(G) + 
  n^{d-1}\cdot\Big(\frac12\,R_0(G) + \frac1{d-1}\,R_1(G)\Big)\\
&+ 
  n^{d-2}\cdot\Big(\frac{d-1}{12}\,R_0(G) + \frac12\,R_1(G)
            + \frac1{d-2}\,R_2(G)\Big)
+ \log n\cdot R_d(G)
  + O(n^{d-3}).
\end{aligned}
$$
\end{proposition}
\begin{remark}
\label{remtriv}
We see that $\tr{}P_n(\log{}B)P_n$
in Theorem~\ref{thm_onedim} and \ref{thm_main} ($G=\log{}B$)
contributes to the leading asymptotic term of order $n^d$, and also to
all lower order terms of order $n^{j}$, $j=d-1,\cdots,1,0,-1,\cdots$,
and to the logarithmic term $\log{n}$, as $n\ra\infty$.
In the classical SSLT the situation is much simpler:
$\log{}B$ is just the Toeplitz matrix of the operator of multiplication
by $\log{}b$, and so $\tr{}P_n(\log{}B)P_n=(2n+1)\widehat{(\log{}b)}_0$.
\end{remark}
Now we are ready to state the two corollaries. 
\begin{corollary}
\label{cor4}
Let $B\in\Psi^0(\S^1)$ have a strictly positive principal
symbol and be such that a certain symbolic norm of $I-B$ is sufficiently small.
Assume also that \eqref{eq_assu} holds.
Then the following holds, as $n\ra\infty$,
$$
 \log\det{}P_nBP_n = c_1\cdot n 
  +  c_{{\rm log}}\cdot\log n
          + c_0 + c_{-1} \cdot \frac1n + O\Big(\frac{1}{n^2}\Big),
$$
where the coefficients are the sums of the corresponding
coefficients from Theorem~\ref{thm_onedim} and Proposition~\ref{prop3}(i).

Assume further that 
$\sigma_0(B)$ 
and $\sub(B)$ do not depend on the direction of $\xi$,
that is $\sigma_0(B)(x,\xi)=b_0(x)$ 
and $\sub(B)(x,\xi)=b_{{\sub}}(x)\,|\xi|^{-1}$, for $(x,\xi)\in{}S^*\S^1$.
Assume also that $b_{-2}=0$. Then the following holds, as $n\ra\infty$,
\begin{equation}
\label{eqmain}
\begin{aligned}
      \log\,&\det{}P_nBP_n = n\cdot
     2 \int_0^{2\pi}\log b_0(x)\,\frac{dx}{2\pi}\\
                       &+ \log n\cdot2\int_0^{2\pi}
       \frac{b_{{\rm sub}}(x)}{b_0(x)}\,\frac{dx}{2\pi}\\ 
                      &+\Bigg(
    \sum_{k=1}^\infty k\,\widehat{(\log b_0)}_k
   \widehat{(\log b_0)}_{-k} 
+ C(\log{}B)+ \gamma\,R_1(\log{}B)+\sum_{l=2}^\infty\zeta(l)\,R_{l}(\log{}B)\Bigg)\\
         &+ \frac1n\cdot \Bigg(    \sum_{k=1}^\infty k\,
  \widehat{(\log b_0)}_k
\big(\widehat{b_{{\rm sub}}/b_0}\big)_{-k}
+\int_0^{2\pi}\bigg[\frac{b_{{\rm sub}}(x)}{b_0(x)}
                           +\bigg(\frac{b_{{\rm sub}}(x)}{b_0(x)}\bigg)^2\,\bigg]
         \,\frac{dx}{2\pi}\Bigg)\\
            &+O\bigg(\frac1{n^2}\bigg),
\end{aligned}
\end{equation}
where $C(\log{}B)$ is given by \eqref{eq_C}.
\end{corollary}
\begin{remark}
In some simple cases, for instance for $b_{{\rm sub}}(x)=\pm\frac12{b_0(x)}$,
the left-hand side in \eqref{eqmain} can be computed explicitly.
The coefficients of $n$, $\log{}n$, and $\frac1n$ on the right in \eqref{eqmain}
in these cases are as expected,
see also Remark~\ref{remmatrix}.
\end{remark}
\begin{corollary}
\label{cor5}
Let $M$ be a Zoll manifold of dimension $d\geq2$.
Assume that $P_n$ and $A$ are as in Theorem~\ref{thm_main}.
Let $B\in\Psi^0(M)$ have a strictly positive principal
symbol and be such that the symbolic norm of $I-B$ is sufficiently small.
For $d=2$, assume in addition \eqref{eq_assu}.
Then the following holds, as $n\ra\infty$,
\begin{equation}
\label{eq_expl}
 \log\det{}P_nBP_n = C_d^{(d)}\cdot n^d  + C_{d-1}^{(d)}\cdot n^{d-1}  
   +C_{d-2}^{(d)}\cdot n^{d-2}  
  +  C_{{\rm log}}^{(d)}\cdot\log n
          + O\big(n^{d-3}\big),
\end{equation}
where the coefficients are the sums of the corresponding
coefficients from Theorem~\ref{thm_main} and Proposition~\ref{prop3}(ii) or (iii).
If one counts the logarithmic term, 
this expansion is fourth order for $d=2,3$
and third order for $d\geq4$. 
\end{corollary}
\begin{remark}
The coefficients $C_d^{(d)}$, $C_{d-1}^{(d)}$, and also $C_{{\rm log}}^{(d)}$, 
$d\in\N$,
have been found in \cite{GOs,GO}.
\end{remark}
\begin{remark}
\label{rem_appl}
The most interesting coefficient in \eqref{eq_expl}
is the constant one, since one can think
of $\exp{}C_0^{(d)}$ as of a regularized determinant of $B\in\Psi^0(M)$,
see \cite{GOs,O2}. The sum 
$$
   \gamma{}R_d(\log{}B)+\sum_{l=2}^\infty\zeta(l)R_{l+d-1}(\log{}B)
$$
will for all $d\in\N$ be a part of $C_0^{(d)}$.
For $d=1$, Corollary~\ref{cor4} gives a full expression for $C_0^{(1)}$.
For $d=2$, Corollary~\ref{cor5} gives a full expression for $C_0^{(2)}$,
which is quite lengthy.
\end{remark}
\begin{remark}
Let us compare the result of Corollary~\ref{cor4} with
a generalization of SSLT to the case of $B$ being an operator
of multiplication by a function $b(x)$ having discontinuities
which is due to H.~Widom and E.~Basor.
In this case $\log{}b(x)$ also has discontinuities, and so the
series $\sum_{k\in\Z}|k|\,|\widehat{(\log{b})}_k|^2$
diverges logarithmically.
The following
third order asymptotic formula holds for 
the operator of multiplication by a piecewise $C^2$
function $b(x)$
\begin{equation}
\label{eqBasor}
       \log\det P_nBP_n = a_1\cdot n+a_2\cdot\log n+a_3+o(1),\qquad n\ra\infty,
\end{equation}
where $a_1$ as in \eqref{eqmain},
the coefficient $a_2$ has been
computed by H.~Widom 
in \cite{W0}, and 
the constant term $a_3$ has been found by E.~Basor in \cite{B}. 
Note that the matrix of $B$ in \eqref{eqBasor} is still Toeplitz,
the logarithmic order of the subleading term being
due to a slower decay of the Fourier coefficients of $b(x)$. 
In our case the matrix of the operator
$B\in\Psi^0(\S^1)$ is {\em not\ }Toeplitz (see Remark~\ref{remmatrix}),
and the logarithm comes from
the contribution of $\sub(B)$.
\end{remark}
\begin{remark}
It would be interesting to find a compact formula
for the constant term in \eqref{eqmain}.
We mention that the constant $a_3$ in \eqref{eqBasor}
found in \cite{B} has a form similar to the one in \eqref{eqmain}.
It contains a ``finite'' term and and an infinite series of
certain integrals multiplied by the values of the Riemann zeta function at 
the points $3,5,\cdots$. Interestingly, an ``invariant'' form
of that series has been found in \cite{W2}. It is written as a single
integral involving the function 
$$
       \Psi(x):=\frac{d}{dx}\,\log \Gamma(x).
$$
This gives the hope that a similar formula can be
found for the constant \eqref{eqmain}.
%
We have done some 
computations trying to find the constant term in \eqref{eqmain},
and the function $\Psi(x)$ has been appearing there.
\end{remark}
\begin{remark} 
\label{remmatrix}
The matrix interpretation of Corollary~\ref{cor4} is as follows.
Assume for simplicity that $B\in\Psi^0(\S^1)$ 
is as in the second part of Corollary~\ref{cor4},
that is $\sigma_0(B)(x,\xi)=b_0(x)$ 
and $\sub(B)(x,\xi)=b_{{\sub}}(x)\,|\xi|^{-1}$, 
for all $(x,\xi)\in{}S^*\S^1$.
Assume also that $b_{-2}=b_{-3}=\cdots=0$.
Let $B_0$ and $B_{{\rm sub}}$ be the operators of
multiplication by $b_0$ and $b_{{\rm sub}}$, respectively.
Let $D$ be the linear operator in $L^2(\S^1)$ such that
$$
        De^{ikx}=\begin{cases}\frac1{|k|}e^{ikx},&|k|\geq1\cr0,&k=0.\end{cases}
$$
Note that this is not a differential, but rather a smoothing operator
of order $-1$.
There is known a correspondence between the classical PsDO's on the
circle and their discrete counterparts, see \cite{M} for details.
By that correspondence, the zeroth order PsDO $B$ we started with,
equals $B_0+B_{{\rm sub}}D$. 
Introduce two Toeplitz matrices,
$\widehat{B}_0:=\{\widehat{(b_0)}_{j-k}\}_{j,k\in\Z}$
and $\widehat{B}_{{\rm sub}}:=
\{\widehat{(b_{{\rm sub}})}_{j-k}\}_{j,k\in\Z}$.
Set also $\widehat{D}:=\diag(\cdots,\frac{1}{3},\frac12,1,0,1,\frac12,\frac{1}{3},\cdots)$.
Then the matrix representation of $B_0+B_{{\rm sub}}D$
is $\widehat{B}_0+\widehat{B}_{{\rm sub}}\cdot{}\widehat{D}$.
Finally, set $\widehat{P}_n=\diag(\cdots,0,1,\cdots,1,0,\cdots)$ ($2n+1$ ones).
We see that Corollary~\ref{cor4} gives a fourth order asymptotics of 
the determinant of the truncated matrix 
$\widehat{P}_n\cdot(\widehat{B}_0+\widehat{B}_{{\rm sub}}\cdot{}\widehat{D})\cdot\widehat{P}_n$.

Now we can reformulate the question 
of finding the constant term in \eqref{eqmain} in purely matrix terms.
Drop the hats and the dots 
for brevity.
Let $C_1$ be a Toeplitz matrix that corresponds to
the operator of multiplication by $b_{{\rm sub}}/b_0$,
and let the matrix $D$ be as above.
Clearly, the matrices $C_1$ and $D$ do not commute.
Assume that the matrix $\log(I-C_1D)$
is well-defined.
The question is to compute the constant coefficient 
in $\tr P_n\log(I-C_1D)P_n$,
or which is the same, the constant coefficient in
\begin{equation}
\label{eqcc}
        \tr{}P_n\log(I-D^{1/2}C_1D^{1/2})P_n, \qquad n\ra\infty.   
\end{equation}
As we have noticed in Remark~\ref{remtriv},
this question is trivial for a Toeplitz matrix $T$ 
in place of $D^{1/2}C_1D^{1/2}$.
\end{remark}
\section{Generalized Hunt--Dyson combinatorial formula}
In this section we state the generalized Hunt--Dyson formula (gHD)
and remind the reader 
the Bohnenblust--Spitzer theorem (BSt).
After that we briefly explain how the gHD
is derived from the BSt 
in \cite{G2}.
We refer to \cite[Section~4]{G2}
for the first step of an independent
proof of the gHD. This first step is a generalization
of F.~J.~Dyson's idea on which the proof of
the usual Hunt--Dyson formula in \cite{K}
is based. It is also explained in \cite[Section~4]{G2}
how to reprove the BSt starting with the gHD.
See \cite[Chapter~2]{G} for the details of the independent
proof of the gHD.

We state the result for the maximum, 
for the corresponding result
for the minimum one should replace the positive parts with the negative
parts.
For $a\in\R$, $n\in\N$ denote
$$
    (a)_+  := \max(0,a),\qquad (a)_+^n:= ( (a)_+)^n.
$$
Fix any $m\in\N$ and assume $a_1,\ldots,a_m\in\R$.
Let $S_m$ be the set of all permutations $\tau$ of
the numbers $1,\ldots,m$. For each
$
     \tau = \genfrac{(}{)}{0pt}{}{1,2,\cdots,m}{\tau_1,\tau_2,\,\cdots,\tau_m}
         \in{}S_m
$
denote
$$
    a_\tau  := (a_{\tau_1},\cdots,a_{\tau_m}).
$$
Introduce the notation
\begin{equation}
\label{CH2_s2e9}
   M_{j}^{}(a_\tau) :=
      \begin{cases}
         \max(0,a_{\tau_1},a_{\tau_1}+a_{\tau_2},\cdots,
        a_{\tau_1}+\cdots+a_{\tau_j}),&j=1,\cdots,m,\cr
         0, &j=0.
      \end{cases}
\end{equation}
Fix any $j=1,\cdots,m$. 
For arbitrary 
$$
           k_1\geq1,\cdots,k_j\geq1,\qquad k_1+\cdots+k_j=m,
$$
we introduce the notation
\begin{equation}
\label{CH2_s1e6prim}
\begin{aligned}
   k_1(a_\tau) &:= a_{\tau_1}+\cdots+a_{\tau_{k_1}}\\
   k_2(a_\tau) &:= a_{\tau_{k_1+1}}+\cdots+
           a_{\tau_{k_1+k_2}}\\
                 &\quad\cdots\\
   k_j(a_\tau) &:= a_{\tau_{k_1+\cdots+k_{j-1}+1}}+\cdots+
         a_{\tau_{k_1+\cdots+k_{j-1}+k_j}}.
\end{aligned}
\end{equation}
Each of $k_l(a_\tau)$, $l=1,\cdots,j$, is a sum
of $k_l$ permuted variables out of $a_{\tau_1},\cdots,a_{\tau_m}$ 
so that each of the permuted variables enters exactly one sum.
Recall that 
$\genfrac{(}{)}{0pt}{}{n}{ l_1,\cdots, l_j}:=\frac{n!}{l_1!\cdots{}l_j!}$
denotes a multinomial coefficient, here $n,j\in\N$, 
$l_1,\cdots,l_j\in\N\cup\{0\}$, and $l_1+\cdots+l_j=n$. 
We are ready to state the generalized Hunt--Dyson formula (gHD).
\begin{theorem}
\label{thm_comb}
For an arbitrary power $n\in\N$,
an arbitrary number of variables $m\in\N$,
and for arbitrary $a_1,\cdots,a_m\in\R$, the following holds 
\begin{equation}
\label{CH2_s2e10}
\begin{aligned}
  \sum_{\tau\in{}S_m}\,
         &\big[ ( M_{m}^{{}}(a_\tau) )^n
                             - ( M_{m-1}^{{}}(a_\tau) )^n \big] \\
         &\!\!\!\!\!\!\!\!\!\!\!\!=  \sum_{\tau\in{}S_m}\, \sum_{j=1}^{\min(m,n)}
         \,\frac1{j!}\,
 \sum_{     \genfrac{}{}{0pt}{}{  k_1,\cdots,k_j \geq 1 }
                  { k_1+\cdots+k_j = m }
          } \,
  \sum_{     \genfrac{}{}{0pt}{}{  l_1,\cdots,l_j \geq 1 }
                  { l_1+\cdots+l_j = n }
           } \,
  \genfrac{(}{)}{0pt}{}{n}{ l_1,\cdots, l_j}\,\frac{(k_1(a_\tau))_+^{l_1}}{k_1}\cdots
                \frac{(k_j(a_\tau))_+^{l_j}}{k_j}.
\end{aligned}
\end{equation}
\end{theorem}
\begin{remark} In the case $n=1$ we obtain
$
           j=1,\quad k_1=m,\quad l_1=1,
$
and \eqref{CH2_s2e10} becomes the usual
Hunt--Dyson combinatorial formula (HD)
\cite[after~(4.8)]{K}
\begin{equation}
\label{CH2_s2e10orig}
\begin{aligned}
  \sum_{\tau\in{}S_m}\,
          \ls \,M_{m}^{{}}(a_\tau) 
                             - M_{m-1}^{{}}(a_\tau)\,\rs 
         &=  \sum_{\tau\in{}S_m}\,
            \frac{(a_{\tau_1}+\cdots+a_{\tau_m})_+}{m}\\
          &=(m-1)!\,(a_1+\cdots+a_m)_+.
\end{aligned}
\end{equation}
\end{remark}
Recall now the statement of the 
Bohnenblust--Spitzer theorem (BSt)
\cite[Theorem~2.2]{S1}. It asserts that for any $m\in\N$ and arbitrary
$a_1,\cdots,a_m\in\R$,
the set $\{M_m(a_{\tau})\}_{\tau\in{}S_m}$
contains the same numbers with the same multiplicities
as the set of sums of positive parts of
the sums of $a_1,\cdots,a_m$, arranged according to the cyclic
representations of all $m!$ permutations. 
This becomes clear if we consider a simple example.
Let us choose $m=3$ and any $a_1,a_2,a_3\in\R$. 
The symmetric group $S_3$
consists of six permutations that
can be written via the cyclic representations as 
$$
    S_3=\big\{(123),\,(132),\,(12)(3),\,(13)(2),\,(23)(1),\,(1)(2)(3)\big\}
$$
In this case the BSt states that the set 
\begin{equation}
\label{eqLHS}
   \big\{\max(0,a_{\sigma_1},a_{\sigma_1}+a_{\sigma_2},
a_{\sigma_1}+a_{\sigma_2}+a_{\sigma_3})\big\}_{\sigma\in{}S_3}
\end{equation}
contains the same numbers with the same multiplicities as the set
\begin{equation}
\label{eqRHS}
\begin{aligned}
    \big\{&(a_1+a_2+a_3)_+,\,(a_1+a_3+a_2)_+,\\
              &(a_1+a_2)_++(a_3)_+,\,(a_1+a_3)_++(a_2)_+,\,
              (a_2+a_3)_++(a_1)_+,\\
              &(a_1)_++(a_2)_++(a_3)_+\big\}.
\end{aligned}
\end{equation}
Note that a certain maximum of zero
and accumulating sums of the permuted variables does not
need to equal the element of the set 
on the right-hand side corresponding to the cyclic representation
of {\em that\ }permutation. 
The statement of the BSt is merely that the whole multisets are identical.
\begin{proof}[Proof of Theorem~\ref{thm_comb}: an outline]
By the BSt, for any $m\in\N$ and arbitrary real $a_1,\cdots,a_m$,
the counterparts of the sets \eqref{eqLHS} and \eqref{eqRHS}
contain the same numbers with the same multiplicities.
Therefore for any function $f:\R_+\ra\C$, the sum 
$\sum_{\tau\in{}S_m} f(M_m^{}(a_\sigma))$
equals the sum of $f$'s values over the counterpart of \eqref{eqRHS}.

It turns out that for the multinomial function $f(t)=t^m$
for any $n\in\N$, a further computation can be performed
which leads to the gHD, see \cite[Section~2]{G2}.

Without going into details let us just mention here that
when we apply the multinomial formula in the sum over 
the analog of \eqref{eqRHS}, some terms
have zero powers. Recall that 
we want {\em all\ }factors to be present in the right-hand side
of the gHD \eqref{CH2_s2e10}. 
However it turns out that the terms
having at least one zero power can be summed together, and their
sum gives {\em exactly\ }the sum of the $n$th power
of the ``previous'' maximum $M_{m-1}(a_\sigma)$ over $S_m$.

At this step we use the principle of inclusion and exclusion 
and the Cauchy and Cayley identities from the theory of partitions. 
\end{proof}
\begin{remark}
The steps of the derivation of the gHD from the BSt can be reversed.
After that having started with the gHD,
we can conclude that for any monomial $f(t)$
its sum over \eqref{eqLHS} equals the sum over \eqref{eqRHS}.
Due to an additional linearity, this actually holds
for an arbitrary {\em polynomial\ }$f(t)$. Now using the polynomial
interpolation we arrive at the BSt, see \cite[Section~4]{G2}.
\end{remark}
\begin{remark}
The independent proof of the gHD in \cite[Chapter~2]{G}
proceeds by induction on the power $n\in\N$.
The base of induction is the usual HD ($n=1$).
In the proof of the inductive step,
the key cancellation of the highest power $n$
of the maximum after taking a sum over $S_m$
follows from a generalization of F.~J.~Dyson's
argument from the proof of the usual Hunt--Dyson formula in \cite{K},
see \cite[Section~4]{G2}.
The proof of the inductive step is however quite technical.
\end{remark}
\begin{remark}
In \cite{RS}, the authors rediscover a version of the usual Hunt--Dyson
formula starting from the much more powerful BSt.
\end{remark}

{\bf Acknowledgments.}
Theorem \ref{thm_main} and Theorem \ref{thm_comb}
appear in the author's Ph.D. thesis. 
I would like to express
a deep gratitude to my thesis adviser Ari Laptev, for suggesting
the problem and his constant attention to the work. 

I would also like to thank Percy Deift for giving a reference \cite{RS}
from which I learned about \cite{S1} and
the Bohnenblust--Spitzer
combinatorial theorem, shortly after having discovered
and proved the generalized Hunt--Dyson formula. 

Finally, I would like to thank the Department of Mathematics
of the Royal Institute of Technology (KTH), Stockholm,
for excellent working conditions 
and a generous financial support
during the whole period of my graduate studies. 
\bibliographystyle{amsalpha}

\end{document}